\documentclass[11pt,a4]{article}
\usepackage{latexsym}
  \usepackage{amssymb}
  \usepackage{times}
\usepackage{graphicx}
\usepackage{hyperref}

  \makeatletter
  \def\section{\@startsection {section}{1}{\z@}{-3.5ex plus -1ex minus
   -.2ex}{2.3ex plus .2ex}{\large\bf}}
  \def\subsection{\@startsection{subsection}{2}{\z@}{-3.25ex plus -1ex minus
   -.2ex}{1.5ex plus .2ex}{\normalsize\bf}}
  \makeatother

  \makeatletter
  \@addtoreset{equation}{section}
  
  \makeatother

  \newtheorem{theorem}{Theorem}[section]
  \newtheorem{proposition}[theorem]{Proposition}

  \newtheorem{definition}[theorem]{Definition}

  \def\qed{%
  \mbox{ }%
  \nolinebreak%
  \hfill%
  \rule{2mm} {2mm}%
  \medbreak%
  \par%
  }
  \renewcommand{\qed}{\hfill$\square$}

\begin{document}
\thispagestyle{empty}
\begin{center}
\null\vspace{-1cm}
\vspace{1.5cm}
{\bf ALMOST CONTACT CURVES IN  TRANS-SASAKIAN 3-MANIFOLDS }\\
\vspace{.5cm}
S. K. Srivastava\footnote{Email: sachink.ddumath@gmail.com}\\
{ Department of Mathematics, Central University of Himachal Pradesh,\\
Dharamshala - 176215,  Himachal Pradesh, India}\\
\end{center}
\vspace{.5cm}
\centerline{\bf Abstract}
This paper is devoted to the study of curvature and torsion of almost contact curves in trans-Sasakian 3-Manifolds. The conditions for the frenet curves to be almost contact curves in trans-Sasakian 3-manifolds have been obtained.
\baselineskip=18pt
\bigskip

\noindent {\bf AMS Subject Classification 2010:} 53D15, 53C25.\\
{\bf Key words:}Almost contact curve, Frenet curve, trans- Sasakian manifold.

\setcounter{footnote}{0}

\section{Introduction}
Almost contact curves play a important role in geometry and topology of contact metric manifolds, e.g. 

\noindent {\it a diffeomorphism of a contact manifold is a contact transformation if and only if any almost contact curves in a domain of it go to almost contact curves} \cite{cb}. Several authors have studied almost contact curves in contact geometry such as \cite{de203, cb, co, cm, jw33, jw54, jiieun, mb}. In \cite{cb} Baikoussis and Blair have studied almost contact curves in contact metric 3-manifold and gave the Frenet 3-frame in this class of manifold. Belkhelfa \emph{et al.} have extended some of the results of \cite{cb} from the Riemannian to the Lorentzian case, and classified all biharmonic almost contac curves in Sasaki-Heisenberg spaces\cite{mb}. A result of Blair [see theorem 8.2,  p.134 \cite{de203}] had been generalized by Welyczko\cite{jw33} to the case of  3-dimensional Quasi-Sasakian manifolds. Moreover the author had also obtained some interesting properties of non-Frenet almost contact curves in normal almost paracontact metric 3-manifolds \cite{jw54}. $\ddot{\rm O}$zg$\ddot{\rm u}$r and Tripathi established necessary and sufficient conditions for almost contact curves having parallel mean curvature vector, proper mean curvature vector, being harmonic and being of type $AW(k), k = 1,2,3$; in $\alpha$- Sasakian manifolds\cite{co}. In \cite{cm} Lee characterized almost contact curves in a Sasakian manifold having the following properties: $(i)$ a pseudo-Hermitian parallel mean curvature vector field $(ii)$ a pseudo-Hermitian proper mean curvature vector field in the normal bundle. Recently, Inoguchi and Lee have studied almost contact curves in normal almost contact metric 3-manifold satisfying $\nabla H = \lambda H \,\,\, {\rm or}\,\,\,\nabla^{\perp} H = \lambda H $  and gave natural equations for {\it planar biminimal} curves \cite{jiieun}.

The purpose of this paper is to investigate the properties of almost contact curves in trans-Sasakian 3-manifolds. This paper is organized as follows: In \S 2 we recall some basic definitions and facts about almost contact metric (in brief a.c.m.) manifolds, trans-Sasakian manifolds and Frenet curves. The curvature, torsion of almost contact curves and the conditions for the frenet curves to be almost contact curves in trans-Sasakian 3-manifolds have been obtained, and finally we construct the examples in \S 3. 
\section{Preliminaries}
\subsection{Contact metric manifolds}
A $(2n+1)$-dimensional differentiable manifold $M$ is said to be an {\it almost contact manifold} if its structural group $GL_{2n + 1}\mathbb{R}$ of linear frame bundle is reducible to ${\rm U}(n)\times\{1\}$  (Gray \cite{jgray}). This is equivalent to existence of a tensor field of type (1, 1), a vector field $\xi$ and a 1-form $\eta$ satisfying
\begin{eqnarray}\label{2.1}
\phi^2 = -I +\eta\otimes\xi, && \eta\left(\xi\right)=1
\end{eqnarray}
From these conditions one can easily obtain
\begin{eqnarray}\label{2.2}
\phi\xi = 0, && \eta o\phi = 0.
\end{eqnarray}
Moreover, scince ${\rm U}(n) \times\{1\}\subset {\rm SO}(2n+1)$, $M$ admits a Riemannian metric $g$ satisfying
\begin{eqnarray}\label{2.3}
g\left(\phi X, \phi Y\right) &=& g(X, Y) - \eta(X)\eta(Y)
\end{eqnarray}
for all $X, Y\in \Gamma (TM).$ Such a metric is called an {\it associated metric} (Sasaki \cite{sasaki}) of the almost contact manifold $M.$ With respect to $g$, $\eta$ is metrically dual to $\xi$, that is
\begin{eqnarray}\label{2.4}
g(X, \xi) &=& \eta (X) 
\end{eqnarray}
A structure $(\phi, \xi, \eta, g)$ on $M$ is called an {\it almost contact metric structure} and a manifold $M$ equipped with an almost contact metric structure is said to be an {\it almost contact metric manifold}. The fundamental 2-form $\Phi$ of the manifold is defined by
\begin{eqnarray}\label{2.5}
\Phi (X, Y) = g(X, \phi Y)
\end{eqnarray}
for all $X, Y \in \Gamma (TM).$ An almost contact metric manifold $M$ is said to be a {\it contact metric manifold} if $\Phi = d\eta .$ Here the exterior derivative $d\eta$ is defined by
\begin{eqnarray}\label{2.6}
d\eta(X, Y) = \frac{1}{2}\left(X\eta(Y) - Y\eta(X) - \eta([X, Y])\right).
\end{eqnarray}
On a contact metric manifold, $\eta$ is {\it contact form}, i.e., $\eta\wedge(d\eta)^n\ne 0$ everywhere on $M$. In particular, $\eta\wedge(d\eta)^n\ne 0$ is a volume element on $M$ so that a contact manifold is orientable. Define a (1, 1) type tensor field $h$ and $l$ by $h= \frac{1}{2}L_{\xi}\phi,\,\,\, lX=R(.,\,\xi)\xi$ , where $L$ denotes the Lie differentiation and $R$ the curvature tensor respectively. The operators $h\,\,{\rm and}\,\, l$ are self-adjoint and satisfy: $h\xi = l\xi= 0$ and $h\phi = -\phi h.$  Also we have $Tr. h = Tr. \phi h = 0.$  Moreover, if $\nabla$ denotes the Levi-Civita connection on $M$, then following formulas holds on a contact metric manifold.
\begin{eqnarray}\label{2.7}
\nabla_{X}\xi = -\phi X - \phi hX.
\end{eqnarray}
\begin{eqnarray}\label{2.71}
l=\phi l\phi -2(h^2 + \phi^2).
\end{eqnarray}
On the direct product manifold $M\times\mathbb{R}$ of an almost contact metric manifold $M$ and the real line $\mathbb{R}$, any tangent vector field can be represented as the form $\left(X, f\frac{d}{dt}\right),$ where $X\in\Gamma(TM)$ and $f$ is a function on $M\times\mathbb{R}$ and $t$ is the cartesian coordinate on the real line $\mathbb{R}.$ \\
Define an almost complex structure $J$ on $M\times\mathbb{R}$ by 
\begin{eqnarray}\label{2.8}
J\left(X, \lambda\frac{d}{dt}\right) = \left(\phi X - \lambda\xi, \eta(X)\frac{d}{dt}\right).
\end{eqnarray}
If $J$ is integrable then $M$ is said to be {\it normal}. Equivalently, $M$ is normal if and only if
\begin{eqnarray}\label{2.9}
[\phi, \phi](X, Y) + 2d\eta(X, Y)\xi = 0,
\end{eqnarray}
where $[\phi, \phi]$ is the Nijenhuis torsion tensor of $\phi$ defined by
\begin{eqnarray}\label{2.10}
[\phi, \phi](X, Y) = [\phi X, \phi Y] + \phi^2[X, Y] - \phi[\phi X, Y] - \phi[X, \phi Y]
\end{eqnarray}
for all $X, Y \in \Gamma (TM).$\\
For an arbitrary almost contact metric 3-manifold $M$, we have (\cite{zo}):
\begin{eqnarray}\label{2.11}
(\nabla_{X}\phi)Y = g(\phi\nabla_{X}\xi, Y)\xi - \eta(Y)\phi\nabla_{X}\xi, X\in\Gamma (TM)
\end{eqnarray}
where $\nabla$ is the Levi-Civita connection on $M$. 

\subsection{ trans-Sasakian manifolds}
This class of manifolds arose in a natural way from the classification of almost contact metric structures and they appear as a natural generalization of both Sasakian and Kenmotsu manifolds. In \cite{ag} Gray Harvella classification of almost Hermite manifolds appear as a class $W_{4}$ of Hermitian manifolds which are closely related to locally conformally K$\ddot {\rm a}$hler manifolds. An almost contact metric structure on a manifold $M \times \mathbb{R}$ belongs to the class $W_{4}$. The class $C_{6}\oplus C_{5}$ \cite{jcm} coincides with the class of trans-Sasakian structure of type $(\alpha, \beta)$.\\
\indent An almost contact metric structure $(\phi, \xi, \eta)$ on a connected manifold $M$ is called {\it trans- Sasakian structure} \cite{jao} if $(M\times\mathbb{R}, J, G)$ belongs to the class $W_{4}$ \cite{ag}, where $J$ is the almost complex structure defined by (\ref{2.8}) and $G$ is the product metric on $M\times\mathbb{R}$. This may be expresses by the condition \cite{debjo}
\begin{eqnarray}\label{2.14}
(\nabla_{X}\phi)Y = \alpha(g(X, Y)\xi - \eta (Y)X) + \beta(g(\phi X, Y)\xi - \eta(Y)\phi X)
\end{eqnarray}
for the smooth functions $\alpha$ and $\beta$ on $M$. Hence we say that the trans - Sasakian structure is of type $(\alpha, \beta)$. From (\ref{2.14}) it follows that 
\begin{eqnarray}\label{2.15}
(\nabla_{X}\eta)(Y) = - \alpha g(\phi X, Y) + \beta g(\phi X, \phi Y)
\end{eqnarray}
We note that trans-Sasakian structure of type $(0, 0),$ $(\alpha, 0)$ and $(0, \beta)$ are the {\it cosymplectic}, $\alpha-$Sasakian and $\beta-$Kenmotsu manifold respectively.
\subsection{Frenet Curves}
Let $( M, g)$ be a Riemannian $n-$manifold with Levi-Civita connection $\nabla$. A unit speed curve $\gamma : I\rightarrow M$ is said to be an{ $r-$\it Frenet curve} \cite{de203} if there exists an orthonormal $r-$frame field $(E_{1} = \gamma ', E_{2},..., E_{r})$ along $\gamma$ such that there exist positive smooth functions $k_{1}, k_{2},..., k_{r-1}$ satisfying 
\begin{eqnarray}
 \nabla_{\gamma '}E_{1} = k_{1}E_{2}, & \nabla_{\gamma '}E_{2} = -k_{1}E_{1} + k_{2}E_{3}, ... , \nabla_{\gamma '}E_{r} = -k_{r-1}E_{r-1}.\nonumber
\end{eqnarray}
The function $k_{r}$ is called the $r-$th curvature of $\gamma$. A Frenet curve is said to be\\
 ${\bullet}$   a {\it geodesic} if $ r = 1,$ i.e., $\nabla_{\gamma '}\gamma ' = 0.$\\
 ${\bullet}$   a {\it Riemannian circle} if $ r = 2,$ and $k_{1}$ is non-zero constant.\\
 ${\bullet}$   a {\it helix} of order $ r $ if $k_{1}, k_{2},..., k_{r-1}$ are constants.\\
 In case $n = 3$, we denote by $(E_{1}, E_{2}, E_{3}) = (T, N, B).$ Then we have the Serret-Frenet equation:
 \begin{eqnarray}
 \nabla_{T}T = kN, & \nabla_{T}N = -kT + \tau B \,\, {\rm{and}} \, \,  \nabla_{T}B = -\tau N
\end{eqnarray}
where $T = \gamma '$.\\
The first curvature $k = k_{1}$ and the second curvature $\tau = k_{2}$ are called the {\it geodesic curvature} and {\it geodesic torsion} of $\gamma$, respectively. The vector field $N$ and $B$ are called the unit normal vector field and binormal vector field of $\gamma$, respectively. 
\section{Almost contact Curves}
Let $\gamma : I\rightarrow M$ be a curve parameterized by arc-length (the {\it natural parametrization}) in an almost contact metric 3-manifold $M$ with Frenet frame $(T, N, B).$
\begin{definition} A Frenet curve $\gamma$ in an almost contact 3-manifold $M$ is said to be an {\it almost contact curve} if it is an integral curve of the contact distribution $D = ker\eta,$ equivalently, $\eta(\gamma ') = 0$.\end{definition}
In particular, when $\eta\wedge(d\eta)^n\ne 0$, almost contact curves are traditionally called {\it Legendre curves} (cf.\cite{cb}).

\noindent We begin with a proposition that will motivate the main result: 
\begin{proposition} Let $M$ be a trans-Sasakian 3-manifold. Then for non-geodesic almost contact curve $\gamma: I\rightarrow M$, curvature ($\kappa$) and torsion ($\tau$) are given by
\begin{eqnarray}\label{4.1}
&& \kappa=\sqrt{\beta^2 + \vartheta^2}
 \end{eqnarray}
\begin{eqnarray}\label{4.2}
&& \tau=\arrowvert{\alpha +\frac{\beta\vartheta' - \beta'\vartheta}{\kappa^2}}\arrowvert
 \end{eqnarray}
\end{proposition}
{\bf Proof:} Let $\gamma$ be an almost contact curve on $M$. Then
\begin{eqnarray}\label{4.3}
&& \nabla_{\gamma'}T = \nabla_{\gamma'}\gamma'=-\beta\xi + \vartheta\phi\gamma'
\end{eqnarray}
for some function $\vartheta.$ The unit normal vector field $N$ is given by
\begin{eqnarray}\label{4.4}
&& N = \frac{1}{\kappa}\nabla_{\gamma'}T = -\frac{\beta}{\kappa}\xi + \frac{\vartheta}{\kappa}\phi\gamma'
\end{eqnarray}
Differentiating (\ref{4.4}) along $\gamma',$ we get
\begin{eqnarray}\label{4.5}
&& \nabla_{\gamma'}N= -\kappa\gamma' + p\xi + q\phi\gamma'
\end{eqnarray}
where
\begin{eqnarray}
p = \frac{\vartheta}{\kappa}\alpha - \frac{\beta'\kappa-\beta\kappa'}{\kappa^2}, && q = \frac{\alpha\beta}{\kappa} + \frac{\vartheta'\kappa-\vartheta\kappa'}{\kappa^2}.
\nonumber
\end{eqnarray}
Here $\beta', \delta'$ and $\kappa'$ are
\begin{eqnarray}
&& \beta'(s) = \frac{d}{ds}\beta(\gamma(s)), \vartheta'(s) = \frac{d}{ds}\vartheta(\gamma(s))\,\,  {\rm and} \,\,\,\kappa'(s) = \frac{d}{ds}\kappa(\gamma(s)) .
\nonumber
\end{eqnarray}
From (\ref{4.4}) and $\tau B = \nabla_{\gamma'}N + \kappa T = p\xi + q\phi\gamma',$ we have (\ref{4.1}) and (\ref{4.2}).\qed
\begin{center} {\bf MAIN RESULT}\end{center}

\begin{theorem}\label{k2}For a Frenet curve $\gamma: I\rightarrow M$ in a trans-Sasakian 3-manifold $M$ with $\alpha \ne 0$ and $\beta \ne 0.$ Set $\sigma = \eta(\gamma ').$ If $\tau = |{l}_{1}\alpha + {l}_{2}\beta + {l}_{3}|$ and at one point of I, $\sigma = \sigma ' = \sigma '' = 0,$ then $\gamma $ is an almost contact curve.\\
Where
$${l}_{1} = \frac{1} {\sqrt {1- \sigma^2}}, \,\,\,{l}_{2} = -\frac{pq\sigma} {\sqrt {1- \sigma^2}(p^2 + q^2)}, \,\,\,{l}_{3} = -\frac{p^2} {(p^2 + q^2)}\gamma\left(\frac{\beta}{p}\right), $$ p is non - zero constant on I and q is certain function on I.
    
\end{theorem}
{\bf Proof:}  Suppose that $\gamma '$ is not collinear with $\xi$ and describe curvature $(\kappa)$ and torsion $(\tau)$ of $\gamma$ on I. We may decompose $\nabla_{\gamma '}\gamma '$ as
\begin{eqnarray}\label{nabla}
\nabla_{\gamma '}\gamma ' = \nabla_{\gamma '} T = \frac{p}{\sqrt {1 - \sigma^2}}\phi\gamma ' + \frac{q}{\sqrt {1 - \sigma^2}}\left(\xi - \sigma\gamma '\right)
\end{eqnarray}
Therefore\\
\begin{eqnarray}\label{k}
k &=& \sqrt{p^2 + q^2}
\end{eqnarray}
is curvature of $\gamma.$\\
Using (\ref{2.11}) and (\ref{nabla}), we have
\begin{eqnarray}\label{sigmadsh}
\sigma ' &=& \gamma ' \left(g\left(\xi, \gamma ' \right)\right) \nonumber\\&=& g\left(\nabla_{\gamma '}\xi, \gamma '\right) + g\left(\xi, \nabla_{\gamma '}\gamma '\right) \nonumber\\  &=& \beta\left(1 - \sigma^2\right) + q\sqrt{1 - \sigma^2}.
\end{eqnarray}
From (\ref{nabla}), we find
\begin{eqnarray}\label{e2}
N &=& \frac{1}{k}\nabla_{\gamma '} T \nonumber\\&=& \frac{p}{k\sqrt {1 - \sigma^2}}\phi\gamma ' + \frac{q}{k\sqrt {1 - \sigma^2}}\left(\xi - \sigma\gamma '\right). 
\end{eqnarray}
Let us write, 
\begin{eqnarray}\label{p1}
&& p_{1} = \frac{p}{k\sqrt {1 - \sigma^2}},  q_{1} = \frac{q}{k\sqrt {1 - \sigma^2}}.
\end{eqnarray}
Then (\ref{e2}) becomes
\begin{eqnarray}\label{e2p1}
N &=& p_{1}\phi\gamma ' + q_{1}\left(\xi - \sigma\gamma '\right).
\end{eqnarray}
From (\ref{k}) and (\ref{sigmadsh}), we compute
\begin{eqnarray}\label{p1dsh}
p_{1} ' = \frac{q\left(p'q - pq'\right)}{k^3\sqrt{1 - \sigma^2}} + \frac{pq\sigma}{k\left(1 - \sigma^2\right)} + \frac{p\sigma\beta}{k\sqrt{1 - \sigma^2}},\nonumber\\
q_{1} ' = \frac{p\left(pq' - p'q\right)}{k^3\sqrt{1 - \sigma^2}} + \frac{q^2\sigma}{k\left(1 - \sigma^2\right)} + \frac{q\sigma\beta}{k\sqrt{1 - \sigma^2}}.
\end{eqnarray}
Differentiating (\ref{e2p1}) along $\gamma '$, we have
\begin{eqnarray}\label{nablae2}
\nabla_{\gamma '}N = p_{1} '\phi\gamma ' + p_{1}\left(\left(\nabla_{\gamma '}\right)\gamma ' + \phi\nabla_{\gamma '}\gamma '\right) + q_{1} '\left(\xi - \sigma\gamma '\right) + q_{1}\left(\nabla_{\gamma '}\xi -  \sigma '\gamma ' - \sigma\nabla_{\gamma '}\gamma '\right).
\end{eqnarray}
Using (\ref{2.14}), (\ref{2.15}), (\ref{nabla}), (\ref{sigmadsh}) and (\ref{4.5}); we get
\begin{eqnarray}
&& \nabla_{\gamma '}N = \left[\frac{q\left(p'q - pq'\right)}{k^3\sqrt{1 - \sigma^2}} + \frac{pq\sigma}{k\left(1 - \sigma^2\right)} + \frac{p\sigma\beta}{k\sqrt{1 - \sigma^2}}q\right]\phi\gamma ' \nonumber\\ && + p_{1}\left[\alpha\left(g\left(\gamma ', \gamma '\right)\xi - \eta\left(\gamma '\right)\gamma '\right) + \beta\left(g\left(\phi\gamma ', \gamma '\right)\xi - \eta\left(\gamma '\right)\phi\gamma '\right) + \frac{p}{\sqrt{1 - \sigma^2}}\left(- \gamma ' + \eta\left(\gamma '\right)\xi\right) - \frac{p\sigma}{\sqrt{1 - \sigma^2}}\phi\gamma '\right]  \nonumber\\ && + \left[\frac{p\left(pq' - p'q\right)}{k^3\sqrt{1 - \sigma^2}} + \frac{q^2\sigma}{k\left(1 - \sigma^2\right)} + \frac{q\sigma\beta}{k\sqrt{1 - \sigma^2}}\right]\left(\xi - \sigma\gamma '\right) \nonumber\\ && + q_{1}\left[-\alpha\phi\gamma ' + \beta\left(\gamma ' - \sigma\xi\right) - \beta\left(1 - \sigma^2\right)\gamma ' - q\sqrt{1 - \sigma^2}\gamma ' - \frac{p\sigma}{\sqrt{1 - \sigma^2}}\phi\gamma ' - \frac{q\sigma}{\sqrt{1 - \sigma^2}}\left(\xi - \sigma\gamma '\right) \right]\nonumber 
\end{eqnarray}
which simplifies to, 
\begin{eqnarray}
&& \nabla_{\gamma '}N =\frac{q}{k\sqrt{1 - \sigma^2}}\left[\frac{\left(p'q - pq'\right)}{k^2} - \frac{p\sigma}{\sqrt{1 - \sigma^2}} - \alpha\right]\phi\gamma ' + \frac{p}{k\sqrt{1 - \sigma^2}}\left[\frac{\left(pq' - p'q\right)}{k^2} + \frac{p\sigma}{\sqrt{1 - \sigma^2}} + \alpha\right]\left(\xi - \sigma\gamma '\right) - k\gamma '
\nonumber 
\end{eqnarray}
or, 
\begin{eqnarray}
&& \nabla_{\gamma '}N =-\frac{q}{k\sqrt{1 - \sigma^2}}\left[\frac{\left(pq' - p'q\right)}{k^2} + \frac{p\sigma}{\sqrt{1 - \sigma^2}} + \alpha\right]\phi\gamma ' + \frac{p}{k\sqrt{1 - \sigma^2}}\left[\frac{\left(pq' - p'q\right)}{k^2} + \frac{p\sigma}{\sqrt{1 - \sigma^2}} + \alpha\right]\left(\xi - \sigma\gamma '\right) - k\gamma '
\nonumber 
\end{eqnarray}
or,
\begin{eqnarray}
 && \nabla_{\gamma '}N + k\gamma ' = \tau B = p_{2}\phi\gamma ' + q_{2}\left(\xi - \sigma\gamma '\right)
\end{eqnarray}
where
\begin{eqnarray}
&& p_{2} = -\frac{q}{k\sqrt{1 - \sigma^2}}\left[\frac{\left(pq' - p'q\right)}{k^2} + \frac{p\sigma}{\sqrt{1 - \sigma^2}} + \alpha\right], \nonumber\\ && q_{2} = \frac{p}{k\sqrt{1 - \sigma^2}}\left[\frac{\left(pq' - p'q\right)}{k^2} + \frac{p\sigma}{\sqrt{1 - \sigma^2}} + \alpha\right]
\end{eqnarray}
and 
\begin{eqnarray}
&& p_{2}^2 + q_{2}^2 = \frac{1}{\left(1 - \sigma^2\right)}\left[\frac{\left(pq' - p'q\right)}{k^2} + \frac{p\sigma}{\sqrt{1 - \sigma^2}} + \alpha\right]^2 \nonumber \\ &&= \tau^2. 
\end{eqnarray}
Since, we have assumed that $\tau = |l_{1}\alpha + l_{2}\beta + l_{3}|.$ Therefore, we have
\begin{eqnarray}
| l_{1}\alpha + l_{2}\beta + l_{3}| = |\frac{1}{\sqrt{1 - \sigma^2}}\left[\frac{\left(pq' - p'q\right)}{k^2} + \frac{p\sigma}{\sqrt{1 - \sigma^2}} + \alpha\right]|
\nonumber
\end{eqnarray}
that is, 
 \begin{eqnarray}\label{sigmadd}
\sigma ''+ \frac{2{\sigma ' }^{2}}{1 - \sigma ^{2}} + p^2\sigma = 0.
\end{eqnarray}
For $\sigma$ not constant, write $\mu = \frac{\sigma'}{p}$. Equation (\ref{sigmadd}) yields
\begin{eqnarray}\label{mudmu}
\mu\frac{d\mu}{d\sigma} + \frac{2\sigma\mu ^2}{1 - \sigma ^2} + \sigma = 0,
\end{eqnarray}
where $p$ is non - zero constant. Integrating (\ref{mudmu}), we have
\begin{eqnarray}\label{mu2}
\mu ^2 = p^2\left(C\left(1 - \sigma^2\right) - 1\right)( 1 - \sigma ^2)
\end{eqnarray}
where $C$ is contant of integration. Using at one point of I, $\sigma = \sigma ' = 0$ and $p\ne 0$ we have $C = 1$.\\
Therefore
\begin{eqnarray}
{\sigma '} ^2  = -p^2\sigma ^2( 1 - \sigma ^2). 
\nonumber
\end{eqnarray}
Recalling $\sigma^2\le 1$, we have $\sigma = 0$, which is a contradiction.\qed

\noindent Let us suppose that $ M=\mathbb R^2\times\mathbb R_{+},$ $\omega:M\rightarrow\mathbb{R}_{+}$
and $(x, y, z)$ be cartesian coordinates in $M$, we define a trans-Sasakian structure on $M$ by
\begin{eqnarray}
&& \xi = \frac{\partial}{\partial z},\,\,\, \eta = dz - ydx\nonumber 
\end{eqnarray}
\[ \phi = \left( \begin{array}{ccc}
0 & -1 & 0 \\
1 & 0 & 0 \\
0 & -y & 0 \end{array} \right),\]
\[ g = \left( \begin{array}{ccc}
\omega +y^2 & 0 & -y\\
0 & \omega & 0\\
-y & 0 & 1 \end{array} \right).\]
Certain almost contact curves in the above class of manifolds are given below:\\
{\bf Example 3.4.} Let us suppose that $\omega=\exp(z),$ then the structure $(\phi, \xi, \eta, g)$ is trans-Sasakian structure of type $\left(\frac{-1}{2\exp (z)}, \frac{1}{2}\right).$ \\{\it A curve $\gamma = \left(\gamma^1,\gamma^2,\gamma^3\right)$ in $M$ is almost contact curve if and only if}
\begin{enumerate}
\item [(i)] $\dot{\gamma}^3 = \gamma^2\,\,\dot{\gamma}^1$
\item [(ii)]  $\left(\dot{\gamma}^1\right)^2 + \left(\dot{\gamma}^2\right)^2= \exp\left(-\gamma^3\right).$
\end{enumerate}
The concrete examples of almost contact curves in $M$  are
\begin{enumerate}
\item [\bf{(3.4.1)}] $\gamma(t) = (1, t, 0), t>0-$ a helix with $\kappa=\tau=1/2.$
\item [\bf{(3.4.2)}] $\gamma(t) = (lnt, 2, 2lnt), t>0-$ a curve with $\kappa=1/2$ and $\tau=1/2t^2.$
\end{enumerate} 
\noindent{\bf Example 3.5.} Let us suppose that $\omega=z,$ then the structure $(\phi, \xi, \eta, g)$ is trans-Sasakian structure of type $\left(\frac{-1}{2z}, \frac{1}{2z}\right).$ \\{\it A curve $\gamma = \left(\gamma^1,\gamma^2,\gamma^3\right)$ in $M$ is almost contact curve if and only if}
\begin{enumerate}
\item [(i)] $\dot{\gamma}^3 = \gamma^2\,\,\dot{\gamma}^1$
\item [(ii)]  $\left(\dot{\gamma}^1\right)^2 + \left(\dot{\gamma}^2\right)^2= \left(\gamma^3\right)^{-1}.$
\end{enumerate}
\noindent The concrete examples of almost contact curves in $M$  are
\begin{enumerate}
\item [\bf{(3.5.1)}] $\gamma(t) = (1, t, 0), t>0-$ a helix with $\kappa=\tau=1/2.$
\item [\bf{(3.5.2)}] $\gamma(t) = (\sqrt{2t}, \sqrt{2t}, t), t>0-$ a generalized helix with $\kappa=\tau=1/2t.$
\end{enumerate} 

\noindent {\bf{Acknowledgment}:} The author wish to express their gratitude to D. E. Blair for helpful comments and remarks in preparing this article.

\end{document}